\documentclass[11pt]{amsproc}%
\usepackage{amssymb}
\usepackage{amsfonts}
\usepackage{amsmath}
\usepackage{graphicx}%
\providecommand{\U}[1]{\protect\rule{.1in}{.1in}}
\theoremstyle{plain}

\newtheorem{corollary}{Corollary}

\newtheorem{lemma}{Lemma}

\newtheorem{remark}{Remark}

\newtheorem{theorem}{Theorem}
\numberwithin{equation}{section}
\begin{document}
\title[Approximation by Lipschitz, $C^{p}$ smooth functions]{Approximation of Lipschitz functions by Lipschitz, $C^{p}$ smooth functions on
weakly compactly generated Banach spaces}
\author{R. Fry}
\curraddr{Department of Mathematics and Statistics, Thompson Rivers University,
Kamloops, BC, CANADA.}
\email{rfry@tru.ca}
\author{L. Keener}
\curraddr{Department of Mathematics and Statistics, University of Northern British
Columbia, Prince George, BC, CANADA.}
\email{keener@unbc.ca}
\thanks{}
\dedicatory{ }
\begin{abstract}
This note corrects a gap and improves results in an earlier paper by the first
named author \cite{F3}.

More precisely, it is shown that on weakly compactly generated Banach spaces
$X$ which admit a $C^{p}$ smooth norm, one can uniformly approximate uniformly
continuous functions $f:X\rightarrow\mathbb{R}$ by Lipschitz, $C^{p}$ smooth
functions$.$ Moreover, there is a constant $C>1$ so that any $\eta$-Lipschitz
function $f:X\rightarrow\mathbb{R}$ can be uniformly approximated by $C\eta
$-Lipschitz, $C^{p}$ smooth functions$.$

This provides a `Lipschitz version' of the classical approximation results of
Godefroy, Troyanski, Whitfield and Zizler.

\end{abstract}
\maketitle

\section{Introduction}

The purpose of this note is to correct a gap in the proof of the main result
of \cite{F3}. Specifically, in the original proof the function $\phi\left(
x\right)  =\left(  Sx,Tx\right)  $ was shown to map from $X$ into the open set
$U\subset l_{\infty}\left(  \mathcal{F}\times\mathbb{N}^{2}\right)  \bigoplus
c_{0}\left(  \mathcal{F}\times\mathbb{N}^{2}\right)  $ as required by Haydon's
theorem (see Theorem 1 below), however; for $F:X\rightarrow\mathbb{R}$
continuous and bounded, and $\left\{  x_{n}^{K}\right\}  \subset X$ dense, it
may be that the map $x\rightarrow\left(  F\left(  x_{n}^{K}\right)  \left(
Sx\right)  _{\left(  K,n,m\right)  },F\left(  x_{n}^{K}\right)  \left(
Tx\right)  _{\left(  K,n,m\right)  }\right)  $ does not, as was thought in
\cite{F3}. The present paper mends this difficulty under the formally stronger
hypothesis that $X$ admit a $C^{p}$ smooth norm rather than merely a
Lipschitz, $C^{p}$ smooth bump function. For weakly compactly generated (WCG)
spaces it is unknown if these two conditions are equivalent. We note, however,
that there are $C_{0}\left(  T\right)  $ spaces, where $T$ is a tree, which
admit $C^{\infty}$ smooth bump functions but no G\^{a}teaux smooth norm (see
\cite{H2}). Despite the formally stronger hypothesis, we have otherwise
improved the results from \cite{F3}. In particular, when the function $f$ to
be approximated has convex domain, we are able to remove the condition imposed
in \cite{F3} that $f$ be bounded, and in addition obtain stronger results when
$f$ is Lipschitz (see Theorem 4 below).

We wish to point out that Theorem 4 of this paper has been used in \cite{HJ1}
to construct $C^{1}$ fine approximations on WCG spaces, and in \cite{S} to
construct smooth extensions of functions from closed subspaces of WCG spaces.
Other reasons for the interest in smooth, Lipschitz approximations includes
their importance in the construction of deleting diffeomorphisms on Banach
spaces (see e.g., \cite{AM}), and their use in smooth variational principles
on infinite dimensional Hilbert manifolds (\cite{AFL}, \cite{AFLR}).

Let us give some background. We are considering the problem of uniformly
approximating continuous, real-valued functions on Banach spaces $X$ by
certain smooth functions. This problem has a long history, beginning with the
work of Kurzweil \cite{K} and continuing through to the present (see e.g.,
\cite{FM}, \cite{DGZ} and further references below). The preferred method for
approaching such smooth approximation problems has been via smooth partitions
of unity. Indeed, the ability to uniformly approximate arbitrary continuous
functions on $X$ by $C^{p}$ smooth functions is equivalent to the existence of
$C^{p}$ smooth partitions of unity on $X$ (see e.g., \cite{DGZ}). In this
vein, for Banach spaces $X$ admitting a $C^{p}$ smooth bump function (a
$C^{p}$ smooth function with bounded, non-empty support), the existence of
$C^{p}$ smooth partitions of unity has been established in fairly wide classes
of spaces. For example, when $X$ is separable this was shown by Bonic and
Frampton \cite{BF}, and this was later generalized to weakly compactly
generated spaces \cite{GTWZ} (see also, \cite{DGZ}, \cite{SS}). Recently in
\cite{HH}, it was shown that in $C\left(  K\right)  $ spaces, for $K$ compact,
the existence of $C^{p}$ smooth bump functions and $C^{p}$ smooth partitions
of unity are equivalent.

One of the drawbacks of employing partitions of unity is that it is very
difficult to arrange for the approximating function to possess nice properties
in addition to basic smoothness. For example, if one wishes the smooth
approximate to be convex or a norm, then other techniques are generally
required (see e.g., \cite{MPVZ}$).$ The situation in which one requires the
approximate to be Lipschitz as well as smooth was addressed in a series of
recent papers, \cite{F1}, \cite{F2}, \cite{AFM}, \cite{AFLR}, \cite{HJ1},
\cite{HJ2}. In particular, it was shown in \cite{AFM} that for separable
Banach spaces $X$ admitting a Lipschitz, $C^{p}$ smooth bump function, the
uniform approximation of uniformly continuous, bounded, real-valued functions
$f$ by Lipschitz, $C^{p}$ smooth functions $g$ is possible (we note that the
uniform continuity of $f$ and the existence of a Lipschitz, $C^{p}$ smooth
bump function on $X$ are necessary). In fact, using a simple argument
motivated by \cite{HJ1} (see also Theorem 4 below), the result from \cite{AFM}
can be used to prove: For $X$ a separable Banach space admitting a Lipschitz,
$C^{p}$ bump, there exists $C>1$ so that for $Y\subset X$ any subset,
$\varepsilon>0,$ and $\eta$-Lipschitz function $f:Y\rightarrow\mathbb{R}$,
there exists a $C\eta$-Lipschitz, $C^{p}$ smooth function $g:X\rightarrow
\mathbb{R}$ with $\left\vert f-g\right\vert <\varepsilon$ on $Y.$ When $Y=X$,
this result has been extended to more general range spaces in \cite{HJ1}. In
\cite{F2} this result is shown to hold for maps having any Banach as range if
we suppose that $X$ has an unconditional Schauder basis$.$ A related result
was shown in \cite{AFLR} where it was proven in particular that for a
separable Hilbert space, the approximate $g$ can be chosen Lipschitz and
$C^{\infty}$ smooth with Lipschitz constant arbitrarily close to the Lipschitz
constant of $f.$

Aside from the $p=1$ case when $X$ is a general Hilbert or superreflexive
space (see e.g., \cite{LL}, \cite{C}), and recent results for $X=c_{0}\left(
\Gamma\right)  $ from \cite{HJ2}, all the results stated above concerning
Lipschitz, $C^{p}$ smooth approximation are for separable Banach spaces$.$ As
indicated above, the purpose of this note is to extend some of these results
to the nonseparable, weakly compactly generated case. In this light, for the
real-valued case, the particular result of \cite{AFM} noted previously can be
seen as a `Lipschitz version' of the classical approximation work of
\cite{BF}, while the present paper can be seen as a `Lipschitz version' of the
(implicit) approximation result of \cite{GTWZ}. In this note we also present a
result on the approximation of Lipschitz functions similar in vein to the
Lipschitz result derived from \cite{AFM} described above. To our knowledge,
even for the $C^{2}$ smooth case in non-separable Hilbert space the results
herein are new. Finally, we show how the result of \cite{GTWZ} can be obtained
from our main result. The entire proof is presented here for the sake of
clarity and completeness.

\section{Main Results}

The notation we use is standard, with $X$ typically denoting a Banach space.
Smoothness here is meant in the Fr\'{e}chet sense and function shall mean
real-valued function. $X$ is said to be \textit{weakly compactly generated
(WCG) }if there exists a weakly compact set $K\subset X$ with $\overline
{\mathbf{span}}\left(  K\right)  =X.$ This class includes the separable and
reflexive Banach spaces.

A Banach space $X$ is said to admit a \textit{separable projectional
resolution of the identity (SPRI)}, if for the first ordinal $\mu$ with
$card\left(  \mu\right)  =dens\left(  X\right)  ,$ there exist continuous
linear projections, $\{Q_{\alpha}:\alpha\in\Gamma\},$ where $\Gamma
=[\omega_{0},\mu],$ so that if we set $R_{\alpha}=\left(  Q_{\alpha+1}%
-Q_{a}\right)  /\left(  \left\Vert Q_{\alpha+1}\right\Vert +\left\Vert
Q_{\alpha}\right\Vert \right)  ,$ we have,

\medskip

(i). $Q_{\alpha}Q_{\beta}=Q_{\min(\alpha,\beta)}$

\medskip

(ii). $(Q_{\alpha+1}-Q_{\alpha})(X)$ is separable for all $\alpha\in\Gamma$

\medskip

(iii). For all $x\in X,\left\{  \left\|  R_{\alpha}(x)\right\|  \right\}
_{\alpha}\in c_{0}(\Gamma)$

\medskip

(iv). For all $x\in X,\;x\in\overline{\text{span}\left\{  R_{\alpha}\left(
x\right)  :\alpha<\mu\right\}  }$

\medskip

\noindent One of the keys to our result is the following fundamental theorem
of Haydon,

\medskip

\begin{theorem}
[H1]For any set $L$ there exists an equivalent norm $\left\Vert \left(
\cdot,\cdot\right)  \right\Vert $ on $l_{\infty}\left(  L\right)  \bigoplus
c_{0}\left(  L\right)  $ such that if $U\left(  L\right)  $ is the open
subset
\[
\left\{  \left(  f,x\right)  \in l_{\infty}\left(  L\right)  \bigoplus
c_{0}\left(  L\right)  :\max\left\{  \left\Vert f\right\Vert _{\infty
},\left\Vert x\right\Vert _{\infty}\right\}  <\left\Vert \left\vert
f\right\vert +\frac{1}{2}\left\vert x\right\vert \right\Vert _{\infty
}\right\}  ,
\]
then $\left\Vert \left(  \cdot,\cdot\right)  \right\Vert $ is $C^{\infty}$
smooth on $U\left(  L\right)  $ and depends locally on only finitely many
non-zero coordinates there.
\end{theorem}

\medskip

\noindent We shall also require the following deep result of Amir and
Lindenstrauss (see also \cite{T}),

\medskip

\begin{theorem}
[AL]If $X$ is a WCG Banach space, then $X$ admits a separable projectional
resolution of the identity.
\end{theorem}

\medskip

\noindent We note that our main result is stated for WCG spaces, but it
applies more generally to any Banach space admitting a separable projectional
resolution of the identity. Such spaces include weakly Lindel\"{o}f determined
spaces, duals of Asplund spaces, and $C\left(  K\right)  $ spaces for $K$ a
Valdivia compact (see \cite{DGZ}). We first establish our result for bounded
functions, then relax this condition later for convex domains.

\begin{theorem}
Let $X$ be a WCG Banach space which admits a $C^{p}$ smooth norm. Let
$\varepsilon>0,$ $G\subset X$ an open subset and $f:G\rightarrow\mathbb{R}$ a
uniformly continuous and bounded function. Then there exists a Lipschitz,
$C^{p}$ smooth function $K$ on $G$ with $\left\vert f\left(  x\right)
-K\left(  x\right)  \right\vert <\varepsilon$ for $x\in G.$
\end{theorem}

\medskip

\textbf{Proof.\ \ }To simplify the proof, we\textbf{\ }shall take $G=X,$
leaving the slight technical adjustments to accommodate a general open subset
$G\subset X$ to the reader. In our use of Theorem 1, we may assume that for
some $A\geq2$ we have $\left\Vert \left(  \cdot,\cdot\right)  \right\Vert
_{\infty}\leq\left\Vert \left(  \cdot,\cdot\right)  \right\Vert \leq
A\left\Vert \left(  \cdot,\cdot\right)  \right\Vert _{\infty},$ where
$\left\Vert \left(  \phi,x\right)  \right\Vert _{\infty}=\max\left\{
\left\Vert \phi\right\Vert _{\infty},\left\Vert x\right\Vert _{\infty
}\right\}  .$ As well, for simplicity we shall assume that $f$ is
$1$-Lipschitz; the case where $f$ is uniformly continuous is similar.

\medskip

\noindent Fix a $C^{p}$ smooth norm $\left\Vert \cdot\right\Vert $ on $X,$ and
let $\varepsilon\in\left(  0,1\right)  .$ Because $f$ is bounded, we may
assume that $3/4\geq f>1/2$ by adding a suitable positive constant and
scaling. Also because $f$ is bounded, it can be uniformly approximated within
$\varepsilon/3A$ by a simple function $\varphi$, and hence for the purposes of
approximation it is enough to work with $\varphi,$ which we do for the
remainder of the proof. We note that by choice of $\varepsilon\ $and $A,\ $we
have $0<\varphi\leq1.$ Let the cardinality of the range of $\varphi$ be $N.$

\medskip

\noindent It will be helpful in the sequel to recall the construction of
$\varphi$ here. We evenly partition $\left[  1/2,3/4\right]  $ into
subintervals $I_{i}=(a_{i},b_{i}]$ with midpoint $m_{i}$ and width
$\Delta=\varepsilon_{1}<\varepsilon/3A.$ Define $\varphi\left(  x\right)
=\sum_{i=1}^{N}m_{i}\chi_{f^{-1}\left(  I_{i}\right)  }.$

\medskip

\noindent We may suppose that $X$ is nonseparable (see the Remark at the end
of this note), and given that $X$ is WCG, by Theorem 2 we let $\left\{
Q_{\alpha}\right\}  _{\alpha\in\Gamma}$ be an SPRI on $X,$ and let
$\mathcal{F}$ be the collection of all finite, non-empty subsets of $\Gamma.$
Our goal shall be to construct appropriate maps $S:X\rightarrow l_{\infty
}\left(  \mathcal{F}\times\mathbb{N}^{2}\right)  $ and $T:X\rightarrow
c_{0}\left(  \mathcal{F}\times\mathbb{N}^{2}\right)  $ for use in applying
Theorem 1.

\medskip

\noindent Using property (ii) of an SPRI, for each $K\in\mathcal{F}$ pick a
dense sequence, $\left\{  x_{n}^{K}\right\}  _{n=1}^{\infty}\subset
X_{K}=\mathbf{span}\left\{  R_{\alpha}\left(  X\right)  :\alpha\in K\right\}
$. From property (iv) of an SPRI, we have that
\[
\mathcal{D}=\left\{  x_{n}^{K}:K\in\mathcal{F},\ n\in\mathbb{N}\right\}
\]
is dense in $X.$

\medskip

\noindent In the folowing lemma, $\varepsilon_{1}$, $\varphi,$ and the $m_{i}$
are defined as given above. For $\delta>0,$ let $\zeta_{\delta}\in C^{\infty
}\left(  \mathbb{R},\left[  0,1\right]  \right)  $ be decreasing and Lipschitz
such that $\zeta_{\delta}\left(  t\right)  =1$ iff $t\leq\delta/32$ and
$\zeta_{\delta}\left(  t\right)  =0$ iff $t\geq\delta/16.$ Using this
notation, we shall require the following technical lemma.

\medskip

\begin{lemma}
For every $x_{n}^{K}\in\mathcal{D}$ there exists an associated $x_{n^{\prime}%
}^{K^{\prime}}\in\mathcal{D}$ and a Lipschitz function $\zeta_{n}^{K}\in
C^{\infty}\left(  \mathbb{R},\left[  0,1\right]  \right)  $ with Lipschitz
constant independent of $\left(  K,n\right)  $ and $\zeta_{n}^{K}\left(
t\right)  =0$ for $t\geq\varepsilon_{1},$ such that for any $y$ with
$\left\Vert y-x_{n}^{K}\right\Vert <\delta/32=\varepsilon_{1}/32^{2}$ (i.e.
$\delta=\varepsilon_{1}/32)$) we have%
\begin{align*}
&  \varphi\left(  x_{n^{\prime}}^{K^{\prime}}\right)  \zeta_{\delta}\left(
\left\Vert y-x_{n}^{K}\right\Vert \right)  \zeta_{n}^{K}\left(
||y-x_{n^{\prime}}^{K^{\prime}}||\right) \\
& \\
&  =\sup\left\{  \varphi\left(  x_{m^{\prime}}^{L^{\prime}}\right)
\zeta_{\delta}\left(  \left\Vert y-x_{m}^{L}\right\Vert \right)  \zeta_{m}%
^{L}\left(  ||y-x_{m^{\prime}}^{L^{\prime}}||\right)  :\left(  L,m\right)
\in\left(  \mathcal{F},\mathbb{N}\right)  \right\}  ,
\end{align*}
where given $\left(  L,m\right)  $ we denote by $\left(  L^{\prime},m^{\prime
}\right)  $ the associated pair.
\end{lemma}

\textbf{Proof.\ \ }Fix $x_{n}^{K}\in\mathcal{D}$ and $y\in X$ with $\left\Vert
y-x_{n}^{K}\right\Vert <\delta/32.$ We assume that $\varphi$ is continuous at
$x_{n}^{K}$. We leave to the reader the simple verifications required when
$\varphi$ is not continuous there. For $i=1,...,N$, define $\rho_{i}=\rho
_{i}\left(  K,n\right)  =\inf\{||x_{n}^{K}-w||:w\in\mathcal{D},\varphi(w)\geq
m_{i}\}$ where $\rho_{i}=\infty$ if there are no such $y$'s. It is clear that
at least one $\rho_{i}$ is finite. We also remark that $\rho_{1}=0,$ as
$\varphi\left(  w\right)  =m_{i}\geq m_{1}$ for some $i.$ Unless otherwise
stated, both $i$ and $\rho_{i}$ are understood to be taken with respect to
$\left(  K,n\right)  $ to ease notation. We shall need the following.

\medskip

\noindent\textbf{Fact.\ }Unless $\rho_{i}=\infty$ or $\rho_{i}=0$, $\rho
_{i+1}\geq\rho_{i}+\varepsilon_{1}$.

\medskip

\noindent\textbf{Proof of Fact.\ }Set $x=x_{n}^{K}$ and suppose that for some
$i$, with $\rho_{i}$ finite and non-zero, $\rho_{i+1}<\rho_{i}+\varepsilon
_{1}.$ Then there is a $z_{1}\in D$ with $\varphi(z_{1})=m_{i+1}$ such that
$||z_{1}-x||=\rho_{i}+\varepsilon_{1}-\eta$ for some $\eta$ satisfying
$0<\eta<\rho_{i}$. Let $z_{2}$ be on the line segment $[x,z_{1}]$ with
$||z_{2}-x||=\rho_{i}-\eta/2.$ There is a $z_{3}\in D$ such that
$||z_{2}-z_{3}||<\eta/2$ and $\varphi(z_{3})\leq m_{i-1}$. Now, since $x$,
$z_{1}$, and $z_{2}$ are collinear,

\medskip%

\begin{align*}
||z_{1}-z_{3}||  &  \leq||(z_{1}-x)-(z_{2}-x)||+||z_{2}-z_{3}||\\
&  =||z_{1}-x||-||z_{2}-x||+||z_{2}-z_{3}||\\
&  <\left(  \rho_{i}+\varepsilon_{1}-\eta\right)  -\left(  \rho_{i}%
-\eta/2\right)  +\eta/2=\varepsilon_{1}\text{.}%
\end{align*}

\medskip

\noindent But $\varphi(z_{1})=\varphi(z_{3})+2\varepsilon_{1},$ by definition
of $m_{i},$ implying $f(z_{1})-f(z_{3})\geq\varepsilon_{1},$ and this violates
the fact that $f$ is $1$-Lipschitz. $\square$

\medskip

\medskip

\noindent Returning to the proof of the lemma, we consider three cases. Recall
that $i=i\left(  K,n\right)  $ and $\rho_{i}=\rho_{i}\left(  K,n\right)  $
unless otherwise stated.

\medskip

\begin{enumerate}
\item If for some $j=j(K,n),$ $\rho_{j(K,n)}\in\left[  \delta,\varepsilon
_{1}\right]  ,$ then define a decreasing, Lipschitz function $\zeta_{n}%
^{K}=\zeta_{1}\in C^{\infty}\left(  \mathbb{R},\left[  0,1\right]  \right)  $
such that $\zeta_{1}\left(  t\right)  =1$ for $t\in\left[  0,\rho_{j}%
-\delta/2\right]  ,$ $\zeta_{1}\left(  t\right)  =0$ for $t\geq\rho_{j}%
-\delta/4.$ Pick $x_{n^{\prime}}^{K^{\prime}}$ with $\left\Vert x_{n}%
^{K}-x_{n^{\prime}}^{K^{\prime}}\right\Vert \in\left[  0,\delta/4\right]  $
and $\varphi\left(  x_{n^{\prime}}^{K^{\prime}}\right)  =m_{j-1}.$ Note that
we may assume Lip$\left(  \zeta_{1}\right)  \leq5/\delta.$

\item If for some $j,$ $\rho_{j}\in\left(  0,\delta\right)  ,$ then define a
decreasing, Lipschitz function $\zeta_{n}^{K}=\zeta_{2}\in C^{\infty}\left(
\mathbb{R},\left[  0,m_{j-1}/m_{j}\right]  \right)  $ such that $\zeta
_{2}\left(  t\right)  =m_{j-1}/m_{j}$ for $t\in\left[  0,\delta+\delta
/4,\right]  ,$ $\zeta_{1}\left(  t\right)  =0$ for $t\geq\delta+\delta/2.$
Observe that since $\rho_{j}\neq0,$ $j>1$ and so $m_{j-1}/m_{j}$ is well
defined. Pick $x_{n^{\prime}}^{K^{\prime}}$ with $\left\Vert x_{n}%
^{K}-x_{n^{\prime}}^{K^{\prime}}\right\Vert \in\left[  \rho_{j},\delta\right]
$ and $\varphi\left(  x_{n^{\prime}}^{K^{\prime}}\right)  =m_{j}.$ Note that
we may assume Lip$\left(  \zeta_{2}\right)  \leq\frac{m_{j-1}}{m_{j}}\frac
{5}{\delta}\leq\frac{m_{1}}{m_{2}}\frac{5}{\delta.}\leq\frac{5}{\delta}.$

\item If there is no $j$ such that $0<\rho_{j}\leq\varepsilon_{1}.$ If
$\rho_{j}=0$ with $j$ the maximal such index and $j>1,$ then take
$x_{n^{\prime}}^{K^{\prime}}=x_{n}^{K}$ and $\zeta_{n}^{K}=\zeta_{2}$. If, on
the other hand, $\rho_{1}=0$ is the largest such index, take $x_{n^{\prime}%
}^{K^{\prime}}=x_{n}^{K}$ and $\zeta_{n}^{K}$ of the form $\zeta_{1}$ using
$\rho_{2}$ in its definition. We handle $\rho_{j}>\varepsilon_{1}$ similarly.
\end{enumerate}

\medskip

\noindent We verify that, with the condition on $y$, we must have%
\begin{align*}
&  \varphi\left(  x_{n^{\prime}}^{K^{\prime}}\right)  \zeta_{\delta}\left(
\left\Vert y-x_{n}^{K}\right\Vert \right)  \zeta_{n}^{K}\left(
||y-x_{n^{\prime}}^{K^{\prime}}||\right) \\
& \\
&  \geq\varphi\left(  x_{m^{\prime}}^{L^{\prime}}\right)  \zeta_{\delta
}\left(  \left\Vert y-x_{m}^{L}\right\Vert \right)  \zeta_{m}^{L}\left(
||y-x_{m^{\prime}}^{L^{\prime}}||\right)
\end{align*}
for any $\left(  L,m\right)  \in\left(  \mathcal{F},\mathbb{N}\right)  $.

\medskip

\noindent The argument depends on the case:

\medskip

\noindent$\left(  \mathbf{1}\right)  \ \ $If case 1 holds for $x_{n}^{K}$,
then $\rho_{j}\left(  K,n\right)  >0,$ $\zeta_{\delta}\left(  \left\Vert
y-x_{n}^{K}\right\Vert \right)  =1$ and $\zeta_{n}^{K}\left(  ||y-x_{n^{\prime
}}^{K^{\prime}}||\right)  =1$ since $||y-x_{n^{\prime}}^{K^{\prime}}%
||\leq\delta/32+\delta/4<\rho_{j}-\delta/2$. Suppose that for some $x_{m}^{L}$,%

\begin{align}
&  \varphi\left(  x_{m^{\prime}}^{L^{\prime}}\right)  \zeta_{\delta}\left(
\left\Vert y-x_{m}^{L}\right\Vert \right)  \zeta_{m}^{L}\left(  \left\Vert
y-x_{m^{\prime}}^{L^{\prime}}\right\Vert \right) \nonumber\\
& \\
&  >\varphi\left(  x_{n^{\prime}}^{K^{\prime}}\right)  \zeta_{\delta}\left(
\left\Vert y-x_{n}^{K}\right\Vert \right)  \zeta_{n}^{K}\left(  \left\Vert
y-x_{m^{\prime}}^{L^{\prime}}\right\Vert \right)  =m_{j-1}\text{.}\nonumber
\end{align}
Then $\left\Vert y-x_{m}^{L}\right\Vert <\delta/16,$ else $\zeta_{\delta
}\left(  \left\Vert y-x_{m}^{L}\right\Vert \right)  =0.$ In particular this
shows that $\left\Vert x_{n}^{K}-x_{m}^{L}\right\Vert <\delta/8.$ Now
$\varphi\left(  x_{m^{\prime}}^{L^{\prime}}\right)  \geq m_{j+1}$ is
untenable. Indeed, suppose this were the case. Since by the Fact above,
$\rho_{j+1}$ cannot be in $[0,\varepsilon_{1})$ it must be that
$||y-x_{m^{\prime}}^{L^{\prime}}||\geq\varepsilon_{1}-\frac{\delta}{32}$. But
then $\zeta_{m}^{L}\left(  ||y-x_{m^{\prime}}^{L^{\prime}}||\right)  =0$
(regardless of the way that $\zeta_{m}^{L}$ is defined). So if $\left(
2.1\right)  $ holds it must be that $\varphi\left(  x_{m^{\prime}}^{L^{\prime
}}\right)  =m_{j}$. We consider the subcases.

\medskip

\begin{itemize}
\item If case 1 holds for $x_{m}^{L}$, then $\rho_{j^{\ast}}\left(
L,m\right)  \in\left[  \delta,\varepsilon_{1}\right]  ,$ where $j^{\ast
}=j^{\ast}\left(  L,m\right)  ,$ and $\varphi\left(  x_{m^{\prime}}%
^{L^{\prime}}\right)  =m_{j\ast-1}$. As $\varphi\left(  x_{m^{\prime}%
}^{L^{\prime}}\right)  =m_{j},$ we must then have $j^{\ast}=j+1.$ Let us show
that this leads to a contradiction. If $j^{\ast}=j+1,$ then by the Fact,
$\rho_{j^{\ast}}(K,n)\geq\rho_{j}(K,n)+\varepsilon_{1}\geq\delta
+\varepsilon_{1}.$ But $\left\Vert x_{n}^{K}-x_{m}^{L}\right\Vert <\delta/8$
easily implies that $\left\vert \rho_{j^{\ast}}(L,m)-\rho_{j^{\ast}%
}(K,n)\right\vert <\delta/8,$ and so $\rho_{j^{\ast}}(L,m)>\frac{7\delta}%
{8}+\varepsilon_{1}\notin\left[  \delta,\varepsilon_{1}\right]  ,$ a contradicition.

\item If case 2 holds for $x_{m}^{L}$, then $m_{j}=\varphi\left(
x_{m^{\prime}}^{L^{\prime}}\right)  =m_{j^{\ast}}\Rightarrow j^{\ast}=j,$ and
since $\zeta_{m}^{L}(t)\leq m_{j^{\ast}-1}/m_{j^{\ast}}=m_{j-1}/m_{j}$ for all
$t$, we must have $\varphi\left(  x_{m^{\prime}}^{L^{\prime}}\right)  \geq
m_{j+1}$ which we have observed is untenable.

\item If case 3 holds for $x_{m}^{L}$ with $\rho_{j^{\ast}}=0$ and $j^{\ast
}>1,$ then $m_{j}=\varphi\left(  x_{m^{\prime}}^{L^{\prime}}\right)
=m_{j^{\ast}}\Rightarrow j=j^{\ast},$ and we argue as in case 2 above. If case
3 holds for $x_{m}^{L}$ where $j^{\ast}=1,$ then $m_{j}=\varphi\left(
x_{m^{\prime}}^{L^{\prime}}\right)  =m_{1}$ implying $j=1$ and so $\rho
_{j}\left(  K,m\right)  =0,$ a contradiction. When $\rho_{j^{\ast}%
}>\varepsilon_{1},$ we proceed similarly.
\end{itemize}

\medskip

\noindent This completes the case 1 analysis.

\medskip

\noindent$\left(  \mathbf{2}\right)  \ \ $If case 2 holds for $x_{n}^{K}$,
then $\zeta_{\delta}\left(  \left\Vert y-x_{n}^{K}\right\Vert \right)  =1$ and
$\zeta_{n}^{K}\left(  ||y-x_{n^{\prime}}^{K^{\prime}}||\right)  =m_{j-1}%
/m_{j}$ since $||y-x_{n^{\prime}}^{K^{\prime}}||\leq\delta/32+\delta
<5\delta/4$. Thus%
\[
\varphi\left(  x_{n^{\prime}}^{K^{\prime}}\right)  \zeta_{\delta}\left(
\left\Vert y-x_{n}^{K}\right\Vert \right)  \zeta_{n}^{K}\left(
||y-x_{n^{\prime}}^{K^{\prime}}||\right)  =m_{j-1}\text{.}%
\]
Suppose that for some $x_{m}^{L}$%
\[
\varphi\left(  x_{m^{\prime}}^{L^{\prime}}\right)  \zeta_{\delta}\left(
\left\Vert y-x_{m}^{L}\right\Vert \right)  \zeta_{m}^{L}\left(
||y-x_{m^{\prime}}^{L^{\prime}}||\right)  >m_{j-1}\text{.}%
\]
Then we must have $\varphi\left(  x_{m^{\prime}}^{L^{\prime}}\right)  =m_{j}$
where again $\varphi\left(  x_{m^{\prime}}^{L^{\prime}}\right)  \geq m_{j+1}$
is untenable.\ From this point the argument is as in case 1.

$\medskip$

\noindent$\left(  \mathbf{3}\right)  \ \ $If $\rho_{j}=0$ with $j>1,$ we argue
as above for case 2. If $\rho_{j}=0$ with $j=1,$ we argue as above for case 1.
The case $\rho_{j}>\varepsilon_{1}$ is similar. $\square$

\medskip

\begin{remark}
\label{remarkvalue}It follows from the proof of Lemma 1, that for any
$x_{n}^{K}$ with associated point $x_{n^{\prime}}^{K^{\prime}}$, and $y\in X$
with $\left\Vert y-x_{n}^{K}\right\Vert <\delta/16,$ we have $\varphi\left(
x_{n^{\prime}}^{K^{\prime}}\right)  \zeta_{n}^{K}\left(  \left\Vert
y-x_{n^{\prime}}^{K^{\prime}}\right\Vert \right)  =m_{j-1},$ for some
$j=j\left(  K,n\right)  .$
\end{remark}

\noindent For the remainder of this note, if $x_{n}^{K}$ is given, we shall
denote the associated point in $\mathcal{D}$ as provided in Lemma 1 by
$x_{n^{\prime}}^{K^{\prime}}.$ We also use the quantities $\varepsilon_{1}$
and $\delta$ as defined above, in the sequel. It is worth noting that
$\zeta_{\delta}\left(  t\right)  \leq1$ and $\zeta_{n}^{K}\leq1$ for all
$\left(  K,n\right)  ;$ facts we shall use later. Let us put
\[
L=\max\left\{  \mathbf{Lip}\left(  \zeta_{\delta}\zeta_{1}\right)
,\mathbf{Lip}\left(  \zeta_{\delta}\zeta_{2}\right)  \right\}  .
\]

\medskip

\noindent Now we define a coordinatewise $C^{p}$ smooth map $S:X\rightarrow
l_{\infty}\left(  \mathcal{F}\times\mathbb{N}^{2}\right)  $ by,

\medskip%

\[
\left(  Sx\right)  _{\left(  K,n,m\right)  }=\varphi\left(  x_{n^{\prime}%
}^{K^{\prime}}\right)  \zeta_{\delta}\left(  \left\Vert x-x_{n}^{K}\right\Vert
\right)  \zeta_{n}^{K}\left(  \left\Vert x-x_{n^{\prime}}^{K^{\prime}%
}\right\Vert \right)
\]

\medskip

\noindent noting that $\left(  Sx\right)  _{\left(  K,n,m\right)  }\leq
\max\varphi\leq1,$ since $\zeta_{\delta}\leq1,\ \zeta_{n}^{K}\leq1.$

\medskip

\noindent Let $\nu\in C^{\infty}\left(  \mathbb{R},\left[  0,1\right]
\right)  $ be such that $\nu\left(  t\right)  =0$ for $t\leq1$, $\nu\left(
t\right)  >0$ for $t>1,$ and $0\leq\nu^{\prime}\leq3.$

\medskip

\noindent Now define $T:X\rightarrow c_{0}\left(  \mathcal{F}\times
\mathbb{N}^{2}\right)  $ by,%

\[
\left(  Tx\right)  _{\left(  K,n,m\right)  }=\left(  \frac{1}{nm\left\vert
K\right\vert }\prod_{\alpha\in K}\nu\left(  m\ \left\Vert R_{\alpha}\left(
x\right)  \right\Vert \right)  \right)  \left(  Sx\right)  _{\left(
K,n,m\right)  }.
\]

\medskip

\noindent Let us first see that $T$ maps into $c_{0}\left(  \mathcal{F}%
\times\mathbb{N}^{2}\right)  .$ Let $\varepsilon^{\prime}>0,$ fix $x\in X$ and
fix $N$ with $1/N<\varepsilon^{\prime}.$ Now $\max\left\{  n,m\right\}  >N$
implies $\left(  Tx\right)  _{\left(  K,n,m\right)  }<\varepsilon^{\prime}$
since $\nu\leq1$ and $\left(  Sx\right)  _{\left(  K,n,m\right)  }\leq1.$
Next, by property (iii) of a SPRI, for each $l\in\mathbb{N},$ let $F_{l}$ be a
finite subset of $\Gamma$ such that $\alpha\notin F_{l}$ implies $\left\Vert
R_{\alpha}\left(  x\right)  \right\Vert <1/l,$ which implies $\nu\left(
l\ \left\Vert R_{\alpha}\left(  x\right)  \right\Vert \right)  =0.$ Finally, let

\medskip%

\[
\mathcal{S}=\left\{  \left(  K,n,m\right)  \in\mathcal{F}\times\mathbb{N}%
^{2}:n,m\leq N,\ K\subset F_{m}\right\}  .
\]

\medskip

\noindent Then $\mathcal{S}$ is finite, and $\left(  K,n,m\right)
\notin\mathcal{S}$ implies $\left(  Tx\right)  _{\left(  K,n,m\right)
}<\varepsilon^{\prime}.$

\medskip

\noindent Next, fix $\left(  K,n,m\right)  $ and consider the coordinate
function $x\rightarrow\left(  Sx\right)  _{\left(  K,n,m\right)  }.$ Observe
that, since $\left\Vert \left\Vert \cdot\right\Vert ^{\prime}\right\Vert
\leq1$ and $\varphi\leq1,$ we have

\medskip%

\begin{align*}
&  \left\Vert \left(  Sx\right)  _{\left(  K,n,m\right)  }^{\prime}\right\Vert
\\
& \\
&  =\left\Vert \varphi\left(  x_{n^{\prime}}^{K^{\prime}}\right)  \left(
\zeta_{\delta}\left(  \left\Vert x-x_{n}^{K}\right\Vert \right)  \zeta_{n}%
^{K}\left(  \left\Vert x-x_{n^{\prime}}^{K^{\prime}}\right\Vert \right)
\right)  ^{\prime}\right\Vert \\
& \\
&  \leq L.
\end{align*}

\medskip

\noindent It follows that for any $x,x^{\prime}\in X,$%

\[
\left\vert \left(  Sx\right)  _{\left(  K,n,m\right)  }-\left(  Sx^{\prime
}\right)  _{\left(  K,n,m\right)  }\right\vert \leq L\left\Vert x-x^{\prime
}\right\Vert ,
\]

\medskip

\noindent and so $S:X\rightarrow l_{\infty}\left(  \mathcal{F}\times
\mathbb{N}^{2}\right)  $ is continuous.

\medskip

\noindent Moreover, since each coordinate function $x\rightarrow\left(
Sx\right)  _{\left(  K,n,m\right)  }$ is Lipschitz with constant independent
of $\left(  K,n,m\right)  ,$ $S:X\rightarrow l_{\infty}\left(  \mathcal{F}%
\times\mathbb{N}^{2}\right)  $ is Lipschitz.

\pagebreak

\noindent Next we have,

\medskip%

\begin{align*}
&  \left(  Tx\right)  _{\left(  K,n,m\right)  }^{\prime}=\frac{1}{nm\left\vert
K\right\vert }\sum_{\beta\in K}\ [\prod_{\alpha\in K\backslash\left\{
\beta\right\}  }\nu\left(  m\ \left\Vert R_{\alpha}\left(  x\right)
\right\Vert \right) \\
& \\
&  \times\ \ \nu^{\prime}\left(  m\ \left\Vert R_{\beta}\left(  x\right)
\right\Vert \right)  \ m\ \left\Vert R_{\beta}\left(  x\right)  \right\Vert
^{\prime}\left(  R_{\beta}^{\prime}\left(  x\right)  \right)  \ ]\left(
Sx\right)  _{\left(  K,n,m\right)  }\\
& \\
&  +\left(  \frac{1}{nm\left\vert K\right\vert }\prod_{\alpha\in K}\nu\left(
m\ \left\Vert R_{\alpha}\left(  x\right)  \right\Vert \right)  \right)
\left(  Sx\right)  _{\left(  K,n,m\right)  }^{\prime},
\end{align*}

\medskip

\noindent from which it follows that%

\begin{align*}
&  \left\Vert \left(  Tx\right)  _{\left(  K,n,m\right)  }^{\prime}\right\Vert
\\
& \\
&  \leq\frac{1}{nm\left\vert K\right\vert }\sum_{\beta\in K}\ [\prod
_{\alpha\neq\beta\in K}\nu\left(  m\ \left\Vert R_{\alpha}\left(  x\right)
\right\Vert \right) \\
& \\
&  \times\ \ \nu^{\prime}\left(  m\ \left\Vert R_{\beta}\left(  x\right)
\right\Vert \right)  \ m\ \left\Vert \left\Vert R_{\beta}\left(  x\right)
\right\Vert ^{\prime}\right\Vert \ \left\Vert R_{\beta}^{\prime}\left(
x\right)  \right\Vert \ ]\left\vert \left(  Sx\right)  _{\left(  K,n,m\right)
}\right\vert \\
& \\
&  +\frac{1}{nm\left\vert K\right\vert }\left(  \prod_{\alpha\in K}\nu\left(
m\ \left\Vert R_{\alpha}\left(  x\right)  \right\Vert \right)  \right)
\left\Vert \left(  Sx\right)  _{\left(  K,n,m\right)  }^{\prime}\right\Vert \\
& \\
&  \leq\frac{1}{nm\left\vert K\right\vert }\sum_{\beta\in K}3m\left\Vert
R_{\beta}\right\Vert +L<3L+L=4L,
\end{align*}

\medskip

\noindent where we have used; $\nu\leq1,$ $\nu^{\prime}\leq3,$ $\left\Vert
\left\Vert \cdot\right\Vert ^{\prime}\right\Vert \leq1,$ and $\left\Vert
R_{\beta}\right\Vert \leq1$ for all $\beta\in\Gamma.$

\medskip

\noindent Hence for each coordinate,
\[
\left\vert \left(  Tx\right)  _{\left(  K,n,m\right)  }-\left(  Tx^{\prime
}\right)  _{\left(  K,n,m\right)  }\right\vert \leq4L\left\Vert x-x^{\prime
}\right\Vert ,
\]
implying $T:X\rightarrow c_{0}\left(  \mathcal{F}\times\mathbb{N}^{2}\right)
$ is continuous and, as above, we have that each coordinate function
$x\rightarrow\left(  Tx\right)  _{\left(  K,n,m\right)  }$ is Lipschitz with
constant independent of $\left(  K,n,m\right)  ,$ and so $T:X\rightarrow
c_{0}\left(  \mathcal{F}\times\mathbb{N}^{2}\right)  $ is also Lipschitz.

\medskip

\noindent Now for each fixed $x\in X$, by property (iv) of an SPRI, there
exists $x_{n}^{K}\in\mathcal{D}$ with $\left\Vert x-x_{n}^{K}\right\Vert
<\delta/32$ and $R_{\alpha}\left(  x\right)  \neq0$ for all $\alpha\in K.$ It
follows from Lemma 1 that for this $K$ and $n,$ and any $m,$
\[
\left(  Sx\right)  _{\left(  K,n,m\right)  }=\varphi\left(  x_{n^{\prime}%
}^{K^{\prime}}\right)  \zeta_{\delta}\left(  \left\Vert x-x_{n}^{K}\right\Vert
\right)  \zeta_{n}^{K}\left(  \left\Vert x-x_{n^{\prime}}^{K^{\prime}%
}\right\Vert \right)  =\left\Vert Sx\right\Vert _{\infty}.
\]

\medskip

\noindent Moreover, from the definition of $\nu,$ for sufficiently large
$m\in\mathbb{N}$ we have $\nu\left(  m\ \left\Vert R_{\alpha}\left(  x\right)
\right\Vert \right)  >0$ for all $\alpha\in K,$ and so for this choice of
$\left(  K,n,m\right)  ,$

\medskip%

\begin{align*}
\left(  Tx\right)  _{\left(  K,n,m\right)  }  &  =\left(  \frac{1}%
{nm\left\vert K\right\vert }\prod_{\alpha\in K}\nu\left(  m\ \left\Vert
R_{\alpha}\left(  x\right)  \right\Vert \right)  \right)  \left(  Sx\right)
_{\left(  K,n,m\right)  }\\
& \\
&  =\left(  \frac{1}{nm\left\vert K\right\vert }\prod_{\alpha\in K}\nu\left(
m\ \left\Vert R_{\alpha}\left(  x\right)  \right\Vert \right)  \right)
\left\Vert Sx\right\Vert _{\infty}>0.
\end{align*}

\medskip

\noindent From the observations immediately above, with the same choice of
$\left(  K,n,m\right)  $ for the given $x,$ it follows that,

\medskip%

\begin{align}
\left\Vert Sx\right\Vert _{\infty}  &  =\left(  Sx\right)  _{\left(
K,n,m\right)  }<\left(  Sx\right)  _{\left(  K,n,m\right)  }+\frac{1}%
{2}\left(  Tx\right)  _{\left(  K,n,m\right)  }\nonumber\\
& \\
&  \leq\left\Vert Sx+\frac{1}{2}Tx\right\Vert _{\infty}.\nonumber
\end{align}

\medskip

\noindent Next, since for any $x\in X$ there exists $\left(  K,n,m\right)  $
with $\left(  Tx\right)  _{\left(  K,n,m\right)  }>0,$ we have for such
$\left(  K,n,m\right)  ,$

\medskip%

\begin{align*}
\left(  Tx\right)  _{\left(  K,n,m\right)  }  &  =\left(  \frac{1}%
{nm\left\vert K\right\vert }\prod_{\alpha\in K}\nu\left(  m\left\Vert
R_{\alpha}\left(  x\right)  \right\Vert \right)  \right)  \left(  Sx\right)
_{\left(  K,n,m\right)  }\leq\left(  Sx\right)  _{\left(  K,n,m\right)  }\\
& \\
&  <\left(  Sx\right)  _{\left(  K,n,m\right)  }+\frac{1}{2}\left(  Tx\right)
_{\left(  K,n,m\right)  },
\end{align*}

\medskip

\noindent and hence $\left(  Tx\right)  _{\left(  K,n,m\right)  }<\left\Vert
Sx+\frac{1}{2}Tx\right\Vert _{\infty}.$ That we may replace $\left(
Tx\right)  _{\left(  K,n,m\right)  }$ in this last inequality with $\left\Vert
Tx\right\Vert _{\infty}$ while maintaining the strictness of the inequality
follows from the fact that $T$ maps into $c_{0}\left(  \mathcal{F}%
\times\mathbb{N}^{2}\right)  ;$ and therefore we have
\begin{equation}
\left\Vert Tx\right\Vert _{\infty}<\left\Vert Sx+\frac{1}{2}Tx\right\Vert
_{\infty}.
\end{equation}

\medskip

\noindent We define $\Phi_{F}:X\rightarrow l_{\infty}\left(  \mathcal{F}%
\times\mathbb{N}^{2}\right)  \bigoplus c_{0}\left(  \mathcal{F}\times
\mathbb{N}^{2}\right)  $ by, $\Phi_{F}\left(  x\right)  =\left(  Sx,Tx\right)
.$ The inequalities $\left(  2.2\right)  \ $and $\left(  2.3\right)  $ show
that $\Phi_{F}$ maps into $U=U\left(  \mathcal{F\times}\mathbb{N}^{2}\right)
$ (see Theorem 1), and moreover, given that both $T$ and $S$ are Lipschitz, we
have that $\Phi_{F}$ is Lipschitz as well.

\medskip

\noindent Next let $\left\Vert \left(  \cdot,\cdot\right)  \right\Vert $ be
the $C^{\infty}$ smooth norm on $U\subset l_{\infty}\left(  \mathcal{F}%
\times\mathbb{N}^{2}\right)  \bigoplus c_{0}\left(  \mathcal{F}\times
\mathbb{N}^{2}\right)  $ as given by Theorem 1. Since as shown $\Phi_{F}$ is
continuous and maps into the open subset $U,$ we have that the composition
$\left\Vert \Phi_{F}\left(  x\right)  \right\Vert $ is $C^{p}$ smooth given
that both $S$ and $T$ are coordinatewise $C^{p}$ smooth, and on $U$ the norm
$\left\Vert \left(  \cdot,\cdot\right)  \right\Vert $ depends locally on only
finitely many non-zero coordinates. We note, $\left\Vert \Phi_{F}\left(
x\right)  \right\Vert \leq A\left\Vert \Phi_{F}\left(  x\right)  \right\Vert
_{\infty}=A\max\left\{  \left\Vert Sx\right\Vert _{\infty},\left\Vert
Tx\right\Vert _{\infty}\right\}  \leq A,$ since $\left\Vert Tx\right\Vert
_{\infty}\leq\left\Vert Sx\right\Vert _{\infty}\leq1.$

\medskip

\noindent Now define $\widehat{S}:X\rightarrow l_{\infty}\left(
\mathcal{F}\times\mathbb{N}^{2}\right)  $ by
\[
\left(  \widehat{S}x\right)  _{\left(  K,n,m\right)  }=\zeta_{\delta}\left(
\left\Vert x-x_{n}^{K}\right\Vert \right)  .
\]
Likewise we define
\[
\left(  \widehat{T}x\right)  _{\left(  K,n,m\right)  }=\left(  \frac
{1}{nm\left\vert K\right\vert }\prod_{\alpha\in K}\nu\left(  m\ \left\Vert
R_{\alpha}\left(  x\right)  \right\Vert \right)  \right)  \left(  \widehat
{S}x\right)  _{\left(  K,n,m\right)  }.
\]
We now analogously define $\Phi:X\rightarrow l_{\infty}\left(  \mathcal{F}%
\times\mathbb{N}^{2}\right)  \bigoplus c_{0}\left(  \mathcal{F}\times
\mathbb{N}^{2}\right)  $ by, $\Phi\left(  x\right)  =\left(  \widehat
{S}x,\widehat{T}x\right)  .$ It is easy to see that again $\Phi$ maps into
$U\left(  \mathcal{F\times}\mathbb{N}^{2}\right)  $ using an argument similar
to that used above for $\Phi_{F}.$ Indeed, for $x$ and $x_{n}^{K}$ with
$\left\Vert x-x_{n}^{K}\right\Vert <\delta/32,$ we have $\left\Vert
\widehat{S}x\right\Vert _{\infty}=\zeta_{\delta}\left(  \left\Vert x-x_{n}%
^{K}\right\Vert \right)  =1.$

\medskip

\noindent Finally define,
\[
K\left(  x\right)  =\frac{\left\Vert \Phi_{F}\left(  x\right)  \right\Vert
}{\left\Vert \Phi\left(  x\right)  \right\Vert }.
\]

\medskip

\noindent As noted above, the numerator of $K$ is $C^{p}$ smooth. Also, for
any $x\in X$ we have $\left\Vert \Phi\left(  x\right)  \right\Vert
\geq\left\Vert \left(  \widehat{S}x,\widehat{T}x\right)  \right\Vert _{\infty
}\geq\left\Vert \widehat{S}x\right\Vert _{\infty}=1,$ and so $K$ is $C^{p}$ smooth.

\medskip

\noindent Now, given that $\Phi_{F}$ is Lipschitz, so is the composition
$\left\Vert \Phi_{F}\left(  x\right)  \right\Vert ,$ and similarly for
$\left\Vert \Phi\left(  x\right)  \right\Vert .$ Since in addition $\left\Vert
\Phi_{F}\left(  x\right)  \right\Vert $ is bounded, and the denominator
$\left\Vert \Phi\left(  x\right)  \right\Vert $ is bounded below by $1,$ it
follows that the quotient function $K$ is Lipschitz.

\medskip

\noindent We finally show that $\left\vert K-\varphi\right\vert <\varepsilon.$
To this end fix $x\in X$ and let,

\medskip%

\[
\mathcal{C}=\left\{  \left(  K,n\right)  \in\mathcal{F\times}\mathbb{N}%
:\left\Vert x-x_{n}^{K}\right\Vert <\delta/16<\varepsilon_{1}<\varepsilon
/3A\right\}  .
\]

\medskip

\noindent Note that if $\left(  K,n\right)  \notin\mathcal{C},$ then $\left(
Sx\right)  _{\left(  K,n,m\right)  }=0$ for all $m,$ from which it follows
that $\Phi_{F}\left(  x\right)  _{\left(  K,n,m\right)  }=\Phi\left(
x\right)  _{\left(  K,n,m\right)  }=0.$

\pagebreak

\noindent Now we estimate (using $\varphi\geq0$),

\medskip%

\begin{align*}
&  \left\vert K\left(  x\right)  -\varphi\left(  x\right)  \right\vert \\
& \\
&  =\left\vert \frac{\left\Vert \Phi_{F}\left(  x\right)  \right\Vert
}{\left\Vert \Phi\left(  x\right)  \right\Vert }-\varphi\left(  x\right)
\frac{\left\Vert \Phi\left(  x\right)  \right\Vert }{\left\Vert \Phi\left(
x\right)  \right\Vert }\right\vert \\
& \\
&  =\left\vert \frac{\left\Vert \Phi_{F}\left(  x\right)  \right\Vert
}{\left\Vert \Phi\left(  x\right)  \right\Vert }-\frac{\left\Vert
\varphi\left(  x\right)  \Phi\left(  x\right)  \right\Vert }{\left\Vert
\Phi\left(  x\right)  \right\Vert }\right\vert \\
& \\
&  \leq\frac{1}{\left\Vert \Phi\left(  x\right)  \right\Vert }\left\Vert
\begin{array}
[c]{c}%
\left(  \left(  \varphi\left(  x_{n^{\prime}}^{K^{\prime}}\right)  \zeta
_{n}^{K}\left(  \left\Vert x-x_{n^{\prime}}^{K^{\prime}}\right\Vert \right)
-\varphi\left(  x\right)  \right)  \zeta_{\delta}\left(  \left\Vert
x-x_{n}^{K}\right\Vert \right)  ,\right. \\
\\
\left.  \left(  \varphi\left(  x_{n^{\prime}}^{K^{\prime}}\right)  \zeta
_{n}^{K}\left(  \left\Vert x-x_{n^{\prime}}^{K^{\prime}}\right\Vert \right)
-\varphi\left(  x\right)  \right)  \left(  \widehat{T}x\right)  _{\left(
K,n,m\right)  }\right)
\end{array}
\right\Vert
\end{align*}

\medskip

\noindent Since only those coordinates in $\mathcal{C}$ survive, we need only
consider $\left(  K,n,m\right)  \in\mathcal{C}\times\mathbb{N}$. Recall that
this implies $\left\Vert x-x_{n}^{K}\right\Vert <\delta/16.$ It will suffice
to consider case 1 and case 2 from Lemm 1; case 3 being similar. Then by
Remark \ref{remarkvalue}, $\varphi\left(  x_{n^{\prime}}^{K^{\prime}}\right)
\zeta_{n}^{K}\left(  \left\Vert x-x_{n^{\prime}}^{K^{\prime}}\right\Vert
\right)  =m_{j-1},$ for some $j=j\left(  K,n\right)  $ and all pairs $\left(
x_{n}^{K},x_{n^{\prime}}^{K^{\prime}}\right)  $ where $\left\Vert x-x_{n}%
^{K}\right\Vert <\delta/16.$ It also follows from the definition of $\rho
_{j}\left(  K,n\right)  $ that $\varphi\left(  x_{n}^{K}\right)  =m_{j-1}.$
Now as $\left\Vert x-x_{n}^{K}\right\Vert <\delta/16,$ from the $1$-Lipschitz
property of $f$ we also have (when defined) that either $\varphi\left(
x\right)  =m_{j-2},$ $m_{j-1},$ or $m_{j}.$ In any event, for $p=0,1,2,$ this
gives%
\[
\left\vert \varphi\left(  x_{n^{\prime}}^{K^{\prime}}\right)  \zeta_{n}%
^{K}\left(  \left\Vert x-x_{n^{\prime}}^{K^{\prime}}\right\Vert \right)
-\varphi\left(  x\right)  \right\vert =\left\vert m_{j-1}-m_{j-p}\right\vert
<\varepsilon_{1}<\varepsilon/3A.
\]
Hence, for $\left(  K,n,m\right)  \in\mathcal{C}\times\mathbb{N}$ we have,%

\begin{align*}
&  \left\vert \varphi\left(  x_{n^{\prime}}^{K^{\prime}}\right)  \zeta_{n}%
^{K}\left(  \left\Vert x-x_{n^{\prime}}^{K^{\prime}}\right\Vert \right)
-\varphi\left(  x\right)  \right\vert \left\vert \zeta_{\delta}\left(
\left\Vert x-x_{n}^{K}\right\Vert \right)  \right\vert \\
& \\
&  \leq\left\vert \varphi\left(  x_{n^{\prime}}^{K^{\prime}}\right)  \zeta
_{n}^{K}\left(  \left\Vert x-x_{n^{\prime}}^{K^{\prime}}\right\Vert \right)
-\varphi\left(  x\right)  \right\vert <\varepsilon/3A.
\end{align*}

\medskip

\noindent Similarly,

\medskip%

\[
\left\vert \left(  \varphi\left(  x_{n^{\prime}}^{K^{\prime}}\right)
\zeta_{n}^{K}\left(  \left\Vert x-x_{n^{\prime}}^{K^{\prime}}\right\Vert
\right)  -\varphi\left(  x\right)  \right)  \left(  \widehat{T}x\right)
_{\left(  K,n,m\right)  }\right\vert <\varepsilon/A.
\]

\medskip

\noindent Now, as $1\leq\left\Vert \Phi\left(  x\right)  \right\Vert $ and
$\left\Vert \cdot\right\Vert \leq A\left\Vert \cdot\right\Vert _{\infty}$, the
estimates above finally gives us $\left\vert K\left(  x\right)  -\varphi
\left(  x\right)  \right\vert <\varepsilon.$ $\blacksquare$

\medskip

\begin{corollary}
[GTWZ]Let $X$ be a WCG Banach space which admits a $C^{p}$ smooth norm. Then
$X$ admits $C^{p}$ smooth partitions of unity.
\end{corollary}

\medskip

\textbf{Proof.\ }Let $A\subset X$ be open and bounded, and set $\rho\left(
x\right)  =\mathbf{dist}\left(  x,X\backslash A\right)  .$ Note that $\rho$ is
uniformly continuous, and since $A$ is bounded, $\rho$ is bounded. Now our
Theorem 3 can be applied to $\rho$ to produce $C^{p}$ smooth, uniform
approximates. Finally, an examination of the proof of Theorem VIII.3.12
\cite{DGZ} shows that the uniform smooth approximation of such $\rho$ is
sufficient to conclude that $X$ admits $C^{p}$ smooth partitions of unity.
$\square$

\medskip

\noindent The next result is based on \cite[Lemma 1]{AFK}, which was motivated
by \cite{HJ1}.

\begin{theorem}
\label{Main}Let $X$ be a WCG Banach space which admits a $C^{p}$ smooth norm.
Then we have:

\begin{enumerate}
\item For every convex subset $Y\subseteq X$, every uniformly continuous
function $f:Y\rightarrow\mathbb{R}$, and every $\varepsilon>0$, there exists a
Lipschitz, $C^{p}$-smooth function $K:X\rightarrow\mathbb{R}$ such that
$\left\vert f\left(  y\right)  -K\left(  y\right)  \right\vert <\varepsilon$
for all $y\in Y.$

\item There exists a constant $C_{0}\geq1$ such that, for every subset
$Y\subseteq X$, every $\eta$-Lipschitz function $f:Y\rightarrow\mathbb{R}$,
and every $\varepsilon>0$, there exists a $C_{0}\eta$-Lipschitz, $C^{p}%
$-smooth function $K:X\rightarrow\mathbb{R}$ such that $\left\vert
f(y)-K(y)\right\vert <\varepsilon$ for all $y\in Y.$
\end{enumerate}
\end{theorem}

\noindent\textbf{Proof.\ } For $\left(  1\right)  ,$ observe that because $f$
is real-valued and $Y$ is convex, by \cite[Proposition 2.2.1 (i)]{BL} $f$ can
be uniformly approximated by a Lipschitz map, and hence it is enough to
establish $\left(  2\right)  .$

\medskip

\noindent The proof of Theorem 3 shows that there is $C\geq1$ such that for
every $1$-Lipschitz function $g:X\rightarrow\lbrack0,10]$, there exists a
$C^{p}$ function $\varphi:X\rightarrow\mathbb{R}$ such that

\begin{enumerate}
\item $|g(x)-\varphi(x)|\leq1/8$ for all $x\in X$

\item $\varphi$ is $C$-Lipschitz.
\end{enumerate}

\noindent Indeed, the proof of Lemma 1 shows that the Lipschitz constants of
the $\zeta_{n}^{K}$ depend only on $\varepsilon$ and the bound of $F,$ from
which it follows that the Lipschitz constant of the smooth approximate $K$ in
Theorem 3 has a likewise dependence.

\noindent We first see that this result remains true for functions $g$ taking
values in $\mathbb{R}$ if we replace $1/8$ with $1$ and we allow $C$ to be
slightly larger. Indeed, by considering the function $h=\theta\circ\varphi$,
where $\theta$ is a $C^{\infty}$ smooth function $\theta:\mathbb{R}%
\rightarrow\lbrack0,10]$ such that $|t-\theta(t)|\leq1/4$ if $t\in
\lbrack0,10]$, $\theta(t)=0$ for $t\leq1/8$, and $\theta(t)=10$ for
$t\geq10-1/8$, we get the following result: there exists $C_{0}:=C\,\text{Lip}%
(\theta)$ such that for every $1$-Lipschitz function $g:X\rightarrow
\lbrack0,10]$ there exists a $C^{p}$ function $h:X\rightarrow\lbrack0,10]$
such that

\begin{enumerate}
\item $|g(x)-h(x)|\leq1/2$ for all $x\in X$

\item $h$ is $C_{0}$-Lipschitz

\item $g(x)=0 \implies h(x)=0$, and $g(y)=10 \implies h(y)=10$.
\end{enumerate}

Now, for a $1$-Lipschitz function $g:X\rightarrow\lbrack0,+\infty)$ we can
write $g(x)=\sum_{n=0}^{\infty}g_{n}(x)$, where
\[
g_{n}(x)=%
\begin{cases}
g(x)-10n & \text{ if }10n\leq g(x)\leq10(n+1),\\
0 & \text{ if }g(x)\leq10n,\\
10 & \text{ if }10(n+1)\leq g(x)
\end{cases}
\]
and the sum is locally finite. The functions $g_{n}$ are clearly $1$-Lipschitz
and take values on the interval $[0,10]$, so there are $C^{p}$ functions
$h_{n}:X\rightarrow\lbrack0,10]$ such that for all $n\in\mathbb{N}$ we have
that $h_{n}$ is $C_{0}$-Lipschitz, $|g_{n}-h_{n}|\leq1/2$, and $h_{n}$ is $0$
or $10$ wherever $g_{n}$ is $0$ or $10$. It is easy to check that the function
$h:X\rightarrow\lbrack0,+\infty)$ defined by $h=\sum_{n=0}^{\infty}h_{n}$ is
$C^{p}$ smooth, $C_{0}$-Lipschitz, and satisfies $|g-h|\leq1$. This argument
shows that there is $C_{0}\geq1$ such that for any $1$-Lipschitz function
$g:X\rightarrow\mathbb{[}0,+\infty)$, there exists a $C^{1}$ function
$h:X\rightarrow\lbrack0,+\infty)$ such that

\begin{enumerate}
\item $|g(x)-h(x)|\leq1$ for all $x\in X$

\item $h$ is $C_{0}$-Lipschitz

\item $g(x)=0 \implies h(x)=0$.
\end{enumerate}

Finally, for an arbitrary $1$-Lipschitz function $g:X\rightarrow\mathbb{R}$,
we can write $g=g^{+}-g^{-}$ and apply this result to find $C^{p}$ smooth,
$C_{0}$-lipschitz functions $h^{+},h^{-}:X\rightarrow\lbrack0,+\infty)$ so
that $h:=h^{+}-h^{-}$ is $C_{0}$-Lipschitz, $C^{p}$ smooth, and satisfies
$|g-h|\leq1$.

Now let us prove the Theorem. By replacing $f$ with the function
\[
x\mapsto\inf_{y\in Y}\{f(y)+\eta\Vert x-y\Vert\},
\]
which is a $\eta$-Lipschitz extension of $f$ to $X$, we may assume that $Y=X$.
Consider the function $g:X\rightarrow\mathbb{R}$ defined by $g(x)=\frac
{1}{\varepsilon}f(\frac{\varepsilon}{\eta}x)$. It is immediately checked that
$g$ is $1$-Lipschitz, so by the result above there exists a $C^{p}$ smooth,
$C_{0}$-Lipschitz function $h$ such that $|g(x)-h(x)|\leq1$ for all $x$, which
implies that the function $K(y):=\varepsilon h(\frac{\eta}{\varepsilon}y)$ is
$C_{0}\eta$-Lipschitz and satisfies $|f(y)-K(y)|\leq\varepsilon$ for all $y\in
X$. $\square$

\medskip

\noindent\textbf{Additional Remarks\ }

\begin{enumerate}
\item The result of \cite{GTWZ} supposes the formally weaker hypothesis that
$X$ admit only a $C^{p}$ smooth bump function rather than a $C^{p}$ smooth
norm. However, as noted earlier, for WCG spaces it is unknown if these
conditions are equivalent.

\item We note that in general it is not possible to establish Theorem 3 from
Corollary 1, since while the existence of $C^{p}$ smooth partitions of unity
implies (in particular) that uniformly continuous functions can be uniformly
approximated by $C^{p}$ smooth functions, these latter functions are not
generally Lipschitz.

\item A version of Theorem 3 for separable $X$ is given in \cite{AFM} (see
also \cite{F1}) which does not use the results of Haydon \cite{H1} or Amir and
Lindenstrauss \cite{AL}.
\end{enumerate}

\medskip

\noindent\textbf{Acknowledgment \ }The authors wish to thank Richard Smith and
Luis S\'{a}nchez Gonz\'{a}lez for reading over the manuscript carefully and
making many suggestions which have improved this note.

\end{document}